




\input amstex
\magnification=\magstep1

\baselineskip=16pt

\documentstyle{amsppt}

\NoBlackBoxes
\settabs 10 \columns

\def\sLo{{\cal L}_{1}}

\def\sLinfty{{\cal L}_{\infty}}
\def\Ball{\rm Ball}

\def\e{\epsilon }

\def\Ball{\hbox{Ball }}
\def \sstwo#1{_{\lower2pt\hbox{${\scriptstyle #1}$}}}
\def\ssone#1{\lower3pt\hbox{${\scriptstyle #1}$}}
\def\ss#1{\lower3pt\hbox{${\scriptstyle #1}$}}
\def\Rn#1{\hbox{{\it I\kern -0.25emR}$\sp{\,{#1}}$}}

\def\qed {\vrule height6pt  width6pt depth0pt}

\def\cal{\Cal}

\def\R{\hbox{{\it I\kern -0.25emR}}}

\def\e{\epsilon }
\def\a{\alpha }
\def\b{\beta }
\def\d{\delta }

\def\gam{\gamma}

\def\longlongrightarrow{\relbar\joinrel\relbar\joinrel
\relbar\joinrel\rightarrow}

\def \fn{\{f_n\}_{n=1}^\infty}

\def\absgn{\{|g_n|\}_{n=1}^\infty}
\def \Gn{\{G_n\}_{n=1}^\infty}

\leftheadtext{ WILLIAM B.~JOHNSON}
\rightheadtext{EXTENSIONS OF $c\ssone{0}$}

\topmatter

\title Extensions of $c\ssone{0}$
\endtitle

\author  William~B.~Johnson 
\endauthor

\thanks
The author was supported in part by NSF DMS 93-06376 and the
Mathematical Sciences Research Institute.
\endthanks

\address W.~B.~Johnson
\newline Department of Mathematics, Texas A\&M University, 
College Station, TX
77843, U.S.A. \endaddress
\email johnson\@math.tamu.edu\endemail

\abstract \  If $X$ is a closed subspace of a Banach space $L$
which embeds into a Banach lattice not containing $\ell_\infty^n$'s
uniformly and $L/X$  contains $\ell_\infty^n$'s
uniformly, then $X$ cannot have local unconditional structure in the
sense of Gordon-Lewis (GL-{\sl l.u.st.}).
\endabstract



\endtopmatter

\vfill\eject

\document

\heading 0. Introduction.\endheading

Fifteen years ago, Bourgain [Bou] gave the first example of an
uncomplemented subspace of an $L_1$ space which is itself isomorphic to
an $L_1$ space. He asked whether there was a ``natural"
example of this phenomenon.  In particular,  if one takes the kernel $X$
of the quotient mapping from $L_1$ onto $c_0$ given by 
$f\mapsto \{\int f\cdot r_n \}_{n=1}^\infty$, where
$\{r_n\}_{n=1}^\infty$ are the Rademacher functions, Bourgain asked
whether $X$ is isomorphic to $L_1$ or whether at least $X$ is a ${\Cal
L}_1$ space.  Of course, this space $X$ is not complemented in $L_1$
because the quotient space $L_1/X \equiv c_0$ does not embed into $L_1$.
Actually, Bourgain attributes these questions to Pisier; at any rate,
both of them as well as e.g. Kisliakov, Zippin, Schechtman, and I
thought about them around that time.  

Recently Kalton and 
Pe\l czy\'nski [KP] solved these problems in the negative.  In fact,
they showed that if $X$ is a subspace of $L_1$ and $c_0$ embeds into
$L_1/X$, then $X$ is uncomplemented in its bidual (so that $X$ is not
isomorphic to an $L_1$ space) and there is an operator from $X$ into a
Hilbert space which is not absolutely summing (so that $X$ is not a 
${\Cal L}_1$ space [LP]).  While lecturing on their results in 1995 and
1996, the authors of [KP] asked whether such an $X$ could have local
unconditional structure ({\sl l.u.st.}), [DPR].  In this note we give a
negative answer to this question and go on to show in Corollary~2.2 that
such an $X$ cannot even have GL-{\sl l.u.st.} [GL].  

The algebraic point of view, important for [KP], is critical for this
paper.  In fact, once one draws the diagram (2.3) and completes it to
(2.4), one realizes that the answer to the Kalton-Pe\l czy\'nski
question is already contained in their paper [KP]! While the proof of
the stronger result that
$X$ fails GL-{\sl l.u.st.} if $L_1/X$ contains a copy of $c_0$ does
use some new analytical lemmas which are generalizations of lemmas in
[KP], no doubt the authors of [KP] would have discovered and
proved them had they looked at (2.4). 

For the most part we use standard Banach space theory terminology, as can
be found in [LT1], [LT2].  However, since the algebraic point of view
is so important for us, in section 1 we introduce some standard
algebraic terminology and rephrase in the language of homological
algebra some known and essentially known results from Banach space
theory. 

I thank Mariusz Wodzicki for reminding me that it is OK to ``think
algebraically", and Alvaro Arias for reading and correcting a
preliminary version of this paper.  

\vfill\eject

\heading 1. Algebraic preliminaries.\endheading

In this section we review some facts  about Banach
spaces used in the sequel, but phrase them in the language of
homological algebra. The analytical facts, except for Proposition~1.7,
which is from [KP], have either been known for  twenty years or are small
generalizations of such facts. The
 algebraic point of view provides a good framework for organizing these
analytical results and makes it much easier to see how to approach the
problem  of Bourgain--Pisier and the related one of 
Kalton-Pe\l czy\'nski mentioned in the introduction.  While the
algebraic  point of view is important in the work of Kalton and 
Pe\l czy\'nski [KP], the language used in [KP] is more standard for
Banach space theory.  Some of what we describe appears in Doma\'nski's
paper [Dom1] and dissertation [Dom2] and the draft of the book of
Castillo and Gonz\'alez [CG].  Also, Kalton himself [Kal]  exposed some
of what we treat in the process of developing a ${{\cal L}_{p}}$-space
theory for
$0<p<1$.

We have included rather more material in this section  than  is
needed for solving the problem of Kalton and Pe\l czy\'nski 
mentioned in the Introduction in the  expectation that the
algebraic point of view will be useful for attacking other
problems in Banach space theory.

The category we work in is usually denoted by $Ban$; the objects are
Banach spaces and the morphisms are bounded linear operators.  A
sequence $\cdots \to X_j \to X_{j+1} \to X_{j+2} \to \cdots$ of
morphisms in
$Ban$ is called {\sl exact\/} provided it is exact in the 
larger Abelian category $Vect$  of vector spaces with linear maps
as morphisms.  This just means the range of each of the bounded
linear operators
$ X_j
\to X_{j+1}$ is the kernel of the suceeding one $X_{j+1} \to
X_{j+2}$.  So  the diagram   
 $0\to X {\buildrel J \over \longrightarrow} L {\buildrel Q\over
\longrightarrow} Y\to 0$
is a {\sl short exact sequence\/} exactly when $J$ is an
isomorphic embedding and
$Q$ is surjective with kernel $JX$; that is, up to the usual 
identifications (here we avoid discussing ``natural isomorphisms"),
$X$ is a subspace of $L$ and $Y$ is the quotient space $L/X$.
The short exact sequence 
$0\to X {\buildrel J \over \longrightarrow} L {\buildrel Q\over
\longrightarrow} Y\to 0$ is said to be an {\sl extension\/} of $Y$
and a
{\sl coextension\/} of $X$. We abuse language by calling both the
sequence itself and  the space
$L$ an extension of $Y$ by $X$ and a coextension of
$X$ by
$Y$.  In Banach space theory it is more common to call $L$ a {\sl
twisted sum\/} of
$X$ with $Y$, but here we shall use the categorical language. 

Given a diagram
$$\matrix
&L & {\buildrel Q\over \longrightarrow}  & Y \cr
&  & & \phantom{u}\uparrow u \cr
& & & Z
\endmatrix\leqno{(1.1)}
$$
we say   $u$ {\sl factors through $Q$\/} or, when $Q$
is understood, {\sl lifts to $L$\/} provided there is an
operator
$\tilde{u}:Z\to L$ making the following diagram commute:
$$\matrix
&L & {\buildrel Q\over \longlongrightarrow}  & Y \phantom{u}\cr
&  & \nwarrow \tilde{u} \phantom{\tilde{u}}&
 \uparrow u
\cr & & &  { Z}\phantom{u}
\endmatrix
$$

Dually, given the diagram 
$$\matrix
 &\phantom{u}X & {\buildrel J\over \longrightarrow}  & L \cr
&u\downarrow &  & \cr
&\phantom{u} Z & & \cr
\endmatrix\leqno{(1.2)}
$$
we say that $u$ {\sl factors through $J$\/} or, when $J$ is
understood,  {\sl extends to $L$\/}, provided there is an operator
$\tilde{u}:L\to Z$ making the following diagram commute:  
$$\matrix
 &\phantom{u}X & {\buildrel J\over \longrightarrow}  & L \cr
&u\downarrow &  \phantom{u}\swarrow \tilde{u}   & \cr
&\phantom{u} Z &  & \cr
\endmatrix
$$
Usually when these concepts are  used,  $Q$ is surjective and $J$
is an isomorphic embedding.

 Let
$\Gn$ be a sequence of finite dimensional spaces which is dense, in
the sense of the Banach-Mazur distance, in the collection of all
finite dimensional spaces, and let $C_p$ be the
$\ell_p$--sum of $\Gn$ when $1\le p \le \infty$ and let $C\ss{0}$
be the
$c\ss{0}$--sum of $\Gn$ (this notation differs slightly from what
we used [J1] when we introduced these spaces).  It is also
convenient to use nonseparable versions of these spaces, so given
an infinite cardinal $\aleph$, let $C_p(\aleph)$ be the
$\ell_p$ sum of $\aleph$ copies of $C_p$ ($c\ss{0} $ sum when
$p=0$).  Actually, separability plays no role in our use of $C_p$
and from a categorical perspective it would be more natural to use
everywhere the $\ell_p$ sum of all finite dimensional subspaces
of
$\ell_\infty$, each repeated $\aleph$-times,  where $\aleph$ is suitably
large, but...

We now come to the definitions of colocal extension and local
lifting which are perhaps not so well known but play an important
role in our investigation  (and, implicitly, in that of [KP]). 
While seemingly particular to the category $Ban$,  Wodzicki
has pointed out that there are analogues of these concepts in some
categories studied by algebraists. 
Referring again to
the lifting diagram $(1.1)$, we say that 
$u$ {\sl locally factors through $Q$\/} or, when $Q$
is understood, {\sl locally lifts to $L$\/}, provided that for
every operator $w:C_1 \to Z$, the composition $uw$ factors
through $Q$.  Notice that this is just an economical  way of saying
that for every finite dimensional subspace $E$ of $Z$, there is a
factorization ${\tilde{u}}\ss{E} $ through $Q$ of the restriction
$u\ss{E}$ of
$u$ to $E$ so that $\sup_E ||{\tilde{u}}\ss{E} ||<\infty$.  

If, in (1.1),  $Q$ is quotient mapping, then every operator from a
$\sLo$ space into $Y$ locally lifts to $L$.

Dually, referring to the diagram $(1.2)$, we say that $u$ {\sl
colocally factors through $J$\/} or, when $J$ is understood,  {\sl
colocally extends to $L$\/}, provided that for every operator
$w:Z\to C_\infty$, the composition $wu$ factors through $J$.  
Notice that this is just an economical  way of saying
that for every finite dimensional quotient space $E$ of $Z$, there
is a factorization ${\widetilde{q\ss{E}u}}  $ through $J$ of the
composition of  $u$ with the quotient mapping 
$q\ss{E}$ of $Z$ onto $E$   so that 
$\sup_E ||{\widetilde{q\ss{E}u}}||<\infty$.

If, in (1.2), $J$ is an isomorphic embedding, then every operator
from $X$ into a $\sLinfty$ space colocally extends to $L$.  While not
obvious from the definition, this follows from Proposition 1.1 and
the fact that the second dual of a $\sLinfty$ space is injective.
In Proposition~1.1 we make use of the fact [FJT] that the identity on any
Banach space colocally factors through the embedding of the space into
its bidual.

Although we do not
really need the spaces $C_1$ and $C_\infty$ in the later sections,
it is interesting to note that they allow the concepts of local lifting
and colocal extensions to be expressed in the language of $Ban$. 
Moreover,   Wodzicki has pointed out that these spaces are
useful in the study of deeper algebraic properties of $Ban$.

From a categorical perspective, the definition of colocal extension
is the ``right" definition since it is evident that it is dual (in
the sense of category theory) to the definition of local lifting.
However, from the perspective of the local theory of Banach
spaces,  probably the most natural definition is item (2) in
Proposition 1.1:

\proclaim {Proposition 1.1} Consider the diagram
$$\matrix
 &\phantom{u}X & {\buildrel J\over \longrightarrow}  & L \cr
&u\downarrow &  & \cr
&\phantom{u} Z & & \cr
\endmatrix\leqno{(1.2)}
$$
The following are equivalent:
\item{(1)} $u$ colocally extends to $L$.
\item{(2)} For every closed subspace $W$ of $L$ containing
$\overline{JX}$ as a finite codimensional subspace, 
there is an operator $u\ss{W}:W\to Z$
so that the   diagram
$$\matrix
 &\phantom{u}X & {\buildrel J\over \longrightarrow}  & W \cr
&u\downarrow & \phantom{uu}\swarrow u\ss{W} & \cr
&\phantom{u} Z & & \cr
\endmatrix
$$
commutes and $\sup_{W} ||u\ss{W}||<\infty$.
\item{(3)} $u^*$ factors through $J^*$.
\item{(4)} $u^{**}$ factors through $J^{**}$.
\item{(5)} The operator 
$X {\buildrel u\over \rightarrow} Z {\buildrel i\over \rightarrow}
Z^{**}$ factors through $J$, where 
$Z {\buildrel i\over \rightarrow} Z^{**}$ is the canonical embedding.
\item{(6)} For every operator $Z {\buildrel w\over \rightarrow} M$
with $M$ a dual space, $wu$ factors through $J$.
\item{(7)} For each cardinal number $\aleph$ and  every operator 
$Z {\buildrel w\over \rightarrow} C_\infty(\aleph)$, $wu$  factors
through $J$.
\endproclaim

\demo{Proof} (3)  implies  (4) by taking adjoints, while the reverse
implication follows from taking adjoints and using the fact that
every dual space is norm one complemented in its bidual.

 (4) $\implies$ (6) follows from the following commutative
diagram (each vertical arrow is the canonical embedding  of the
space  into its bidual) and the complementation of $M$ in
$M^{**}$:
$$
\matrix
&X&{\buildrel u\over \longrightarrow}&Z&
{\buildrel w\over \longrightarrow}&M\cr
&\downarrow&&\downarrow&&\downarrow\cr
&\phantom{**}X^{**}&{\buildrel u^{**}\over
\longrightarrow}&Z^{**}& {\buildrel w^{**}\over
\longrightarrow}&M^{**}\cr
&J^{**}\downarrow \phantom{J^{**}}&&\uparrow\cr
&\phantom{**}L^{**}&=&\phantom{**}L^{**}\cr
\endmatrix
$$

(6) implies (7) because $C_\infty(\aleph)=C_1(\aleph)^*$, while
(6) $\implies$ (5) is formal.

(7) $\implies$ (1) is formal and (1) $\implies$ (7) is essentially
obvious.

The implication (2)  $\implies$ (3) uses the ``Lindenstrauss
compactness method" and involves only a small variation of an argument
in [J2], so we just outline the proof.  Extend each of the operators 
$u\ss{W}$ to (nonlinear, discontinous) mappings $v\ss{W}$ from $L$ to
$Z$ by defining
$v\ss{W}(y)$ to be $0$ when $y$ is not in $W$.  The $W$'s are
directed by inclusion and thereby generate a net of functions from
$Z^*$ into $\Rn{L}$ defined by $v^{\#}_W(z^*)(y)=z^*(v\ss{W}y)$.  It
is easy to verify that the net $\{v^{\#}_W\}$ has a cluster point
$v:Z^*\to \Rn{L}$ in the product space $\left(\Rn{L}\right)^{Z^*}$
and that $v$ is in fact a bounded linear operator from $Z^* $ into
$L^*$, and that $u^*=J^*v$.

The aforementioned Lindenstrauss compactness method was used in the
same way in [J2] to prove that if $M$ is a separable Banach space, then
$M^*$ is isometrically isomorphic to a norm one complemented subspace
of
$C_\infty$. (Actually, that is why we outlined the proof of
(2)  $\implies$ (3) rather than the slightly simpler direct proof of
(2)  $\implies$ (4).)   A similar argument yields that if the density
character of
$M$ is $\aleph$, then
$M^*$ is isometrically isomorphic to a norm one complemented subspace
of
$C_\infty(\aleph)$.  This gives (7) $\implies$ (6).

The implication (5) $\implies $ (2) follows from Lemma~2.9 in [FJT]
(which is, in turn, a simple consequence of the Principle of Local
Reflexivity [LR] in the form given in [JRZ]),
which says that condition (2) is true in the special case when
$L$ is the space $X^{**}$, $J$ is the canonical embedding,  and $u$
is the identity operator on $X$.  Indeed, let $\a:L\to Z^{**}$ 
satisfy the factorization identity $\a J = iu$, and for $W$ as in
item (2), let $Z_{W}$ be the linear span in $Z^{**}$ of 
$(iZ)\cup (\a W)$, and notice that $iZ$ has finite codimension in
$Z_W$ because $\overline{JX}$ has finite codimension in $W$.  Given
$\e>0$, Lemma~2.9 in [FJT] says that there is an operator 
$P_W:Z_W\to Z$ so that $\left(P_{W}\right) i=I_Z$
and $||P_W||< 3+\e$.  Setting $u\ss{W}=(P_W)\a\ss{|W}$, we see that
$ u=u\ss{W} J$ and $\sup_{W} ||u\ss{W}||< (3+\e)||\a||$.
\hfill\qed\enddemo

In order to characterize when the operator $u$ in (1.1)
locally lifts to $L$, it is convenient to introduce a weaker
concept of factorization.  In (1.1), say that $u$ {\sl
approximately factors through $Q$\/} or, when $Q$ is
understood, {\sl approximately lifts to $L$}, provided that
for each $\e>0$ there is an operator $u\ss{\e}:Z\to L$ so
that $||u-Qu\ss{\e}||<\e$ and $\sup\ss{\e}
||u\ss{\e}||<\infty$.  Similarly, say that $u$ {\sl
approximately locally factors through $Q$\/} or, when $Q$ is
understood, {\sl approximately locally lifts to $L$},
provided that for every operator $w:C_1\to Z$, the composition
$uw$ approximately factors through $Q$.  This is equivalent
to saying that for every finite dimenional subspace $E$ and
$\e>0$, there is an operator $u\ss{E,\e}:Z\to L$ so that 
$||u-Qu\ss{E,\e}||<\e$ and 
$\sup\ss{E,\e} ||u\ss{E,\e}||<\infty$.  

For  a typical example of an operator which approximately
factors but does not factor, set in (1.1) $L=\ell_1$,
$Y=\ell_2$, $Z=\R$,  let $Q$ be the linear extension of
the the operator which takes the $n$-th unit basis vector
$e_n$ in $\ell_1$ to $e_1+{1\over n} e_{n+1}$ in $\ell_2$,
and define $u(t)=te_1$.

\proclaim{Proposition 1.2} Consider the diagram
$$\matrix
&L & {\buildrel Q\over \longrightarrow}  & Y \cr
&  & & \phantom{u}\uparrow u \cr
& & & Z
\endmatrix\leqno{(1.1)}
$$
The following are equivalent:
\item{(1)} $u$  approximately locally lifts to $L$.
\item{(2)} For every quotient 
$L{\buildrel {q\ss{W}} \over \longrightarrow} W$ of
$L$ by a finite codimensional subspace of $\ker Q$ and every
$\e>0$,  there is an operator $u\ss{W,\e}:Z\to W$
so that 
$||Q_Wu\ss{W,\e}-u||<\e$
 and $\sup_{W,\e} ||u\ss{W,\e}||<\infty$, where $W{\buildrel
Q_W\over \longrightarrow} Y$ is the mapping induced by $Q$.
\item{(3)} $u^*$ factors through $Q^*$.
\item{(4)} $u^{**}$ factors through $Q^{**}$.
\item{(5)} The operator 
$iu$ factors through $Q^{**}$, where 
$Y {\buildrel i\over \rightarrow} Y^{**}$ is the canonical embedding.
\endproclaim

\demo{Proof}  As in the proof of Proposition~1.1, (3) and (4) are
easily seen to be equivalent, and (1) $\implies$ (3) (or (1)
$\implies$ (4)) follows from a simple compactness argument. (4)
$\implies$ (5) is formal.

For  (5) $\implies$ (1), get 
$Z{\buildrel w\over \longrightarrow} L^{**}$ so that $Q^{**} w = i
u$ and fix a finite dimensional subspace $E$ of $Z$. By the
principle of local reflexivity, there is a net $\{v\ss{\d}\}$ of
operators from
$wE$ into $L$ so that $\lim\ss{\d} ||v\ss{\d}|| =1$
and 
$\{v\ss{\d} wz\}$ weak$^*$ converges to $wz$ for each $z$ in
$E$.  Since $Q^{**}$ is weak$^*$ continuous and extends $Q$, 
$\{iQv\ss{\d} wz\}$ weak$^*$ converges in $Y^{**}$ to $Q^{**}wz$
for each $z$ in $E$. But for $z$ in $E$, 
$Q^{**}wz=iuz$,  so in fact 
$\{Qv\ss{\d} wz\}$ converges weakly in $Y$ to $uz$.  Therefore we
can get a net of far out convex combinations of $\{v\ss{\d}\}$,
which we continue to denote by $\{v\ss{\d}\}$, so that for each
$z$ in $E$, 
$\lim\ss{\d} ||Q v\ss{\d}wz - uz||=0$, and hence even
$\lim\ss{\d} ||Q v\ss{\d}w\ss{|E} - u\ss{|E}||=0$.  

We included  (2) mostly because it is the  ``approximate" dual
condition to item (2) in Proposition~1.1 and so omit the proof
that it is equivalent to the other conditions in Proposition~1.2.
\hfill\qed\enddemo

Proposition~1.2 combines with the Proposition~1.3 to give a
characterization of when in (1.1) $u$ locally factors through $Q$
when $Q$ has closed range (the case of interest to us in the next
section).

\proclaim{Proposition 1.3} Consider the diagram:
$$\matrix
&L & {\buildrel Q\over \longrightarrow}  & Y \cr
&  & & \phantom{u}\uparrow u \cr
& & & Z
\endmatrix\leqno{(1.1)}
$$
If $Q$ has closed range and $u$ approximately locally factors through
$Q$, then $u$ locally factors through $Q$.
\endproclaim

\demo{Proof} It is clear from the definition that if an operator
$w$ approximately locally factors through an operator $v$, then the
range of $w$ is contained in the closure of the range of $v$. 
Consequently, since $Q$ has closed range, we can assume that $Q$ is
surjective.

Note that there is a constant $C_n$ so that for every
$n$-dimensional subspace $F$ of $Y$, there is an operator
$w\ss{F}:F\to L$ so that $||w\ss{F}||\le C_n$ and $Qw\ss{F}=I_F$. 
Indeed, since $Q$ is surjective, $Q\Ball(L)\supset \d \Ball(Y)$ for
some $\d>0$.  Take in $F$ an Auerbach basis
$\{y_{j},y^*_j\}_{j=1}^n$; that is,
$y^*_{j}(y_{i})=\d\ss{i,j}$ and $||y^*_{j}||=1=||y_{j}||$; and
choose $x_j$ in $L$ with $||x_j||\le {1\over \d}$ and $Qx_j=y_j$.
Set $w\ss{F}y_j=x_j$ and extend linearly to $F$.  Then
$||w\ss{F}||\le{n\over \d} $ and $Qw\ss{F}=I_F$. 

Choose $C$ so that for every finite dimensional subspace $E$ of $Z$
and $\e>0$, there is an operator $u\ss{E,\e}:E\to L$ so that
$||Qu\ss{E,\e}-u\ss{|E}||<\e$ and $||u\ss{E,\e}||\le C$.  Fix a
finite dimensional subspace $E$ of $Z$ and set $n=\dim E$.  Given
$\e>0$, set $F_\e=\left(Qu\ss{E,\e}-u\right)E$.  Then $v\equiv
u\ss{E,\e}-w\ss{F_\e}\left(Qu\ss{E,\e}-u\right):E\to L$  satisfies
$||v||\le C+\e C_n$ and $Qv=u\ss{|E}$\,.  \hfill\qed\enddemo

A short exact sequence 
$0\to X {\buildrel J \over \longrightarrow} L {\buildrel Q\over
\longrightarrow} Y\to 0$ is said to {\sl split\/} provided the
identity on $Y$ lifts to
$L$.  (Sometimes we abuse language by saying that the extension,
$L$, of $Y$ by $X$  splits.)  This is equivalent to saying that the
identity on
$X$ extends to $L$ which is just to say that $JX$ is a complemented
subspace of $L$.  Say that 
$0\to X {\buildrel J \over \longrightarrow} L {\buildrel Q\over
\longrightarrow} Y\to 0$ {\sl locally splits\/} provided that the
identity on
$Y$ locally lifts to $L$.  The following Corollary, which is an
immediate consequence of  Proposition~1.1, Proposition~1.2, and
Proposition~1.3, gives several equivalents to the concept of local
splitting.  Whatever novelty there may be in Propositions~1.2--1.3,
most of Corollary~1.4 is in the literature.  In particular,
$(1)\implies (3)$ is Proposition~1 in [J2]. That $(4)$ implies the
version of $(2)$ in Proposition~1.1(2) is (as noted in the proof of
Proposition~1.1) essentially Lemma~2.9 in [FJT]; moreover, the
equivalence of $(4)$ with 1.1(2) is part of Theorem~3.5 in [Kal].

\proclaim{Corollary 1.4} The following are equivalent for the
short exact sequence \break
$0\to X {\buildrel J \over \longrightarrow} L {\buildrel Q\over
\longrightarrow} Y\to 0$:
\item{(1)} The sequence locally splits.
\item{(2)} The identity on $X$ colocally extends to $L$.
\item{(3)} The short exact sequence
$0\leftarrow X^* {\buildrel J^* \over \longleftarrow} L^* {\buildrel
Q^*\over
\longleftarrow} Y^*\leftarrow 0$ splits.
\item{(4)} The short exact sequence
$0\to X^{**} {\buildrel J^{**} \over \longrightarrow} L^{**}
{\buildrel Q^{**} \over \longrightarrow} Y^{**}\to 0$ splits.
\endproclaim

The last categorical concepts we mention are those of
pushouts and pullbacks.  Given the diagram (1.1), a {\sl
pullback\/} of it is a commutative diagram
$$\matrix
&L & {\buildrel Q\over \longrightarrow}  & Y \cr
&\a \uparrow\phantom{\a} & & \phantom{u}\uparrow u \cr
& W &{\buildrel \b \over\longrightarrow} & Z
\endmatrix\leqno{(1.3)}
$$
which satisfies the minimality condition that if 
$$\matrix
&L & {\buildrel Q\over \longrightarrow}  & Y \cr
&\a_1 \uparrow\phantom{\a} & & \phantom{u}\uparrow u \cr
& W_1 &{\buildrel \b_1 \over\longrightarrow} & Z
\endmatrix
$$
is another commutative diagram, then there is a unique morphism 
$W_1 {\buildrel w \over\longrightarrow} W$ so that $\a_1=\a w$
and $\b_1=\b w$.  In any category pullbacks are unique in an
obvious sense whenever they exist.  Pullbacks of course do exist
in $Ban$: Given (1.1),   $W$ in (1.3) is the subspace of 
$L\oplus_\infty Z$  of all pairs $(x,z)$ for which 
$Qx=uz$. The operator $\a$ (respectively, $\b$) is the restriction to
$W$ of the coordinate projection from  $L\oplus_\infty Z$ onto $L$
(respectively, $Z$). We call this the {\sl canonical
pullback construction.\/} Forgetting norms, this is same
construction that is used to build pullbacks in the Abelian
category $Vect$, so general categorical principles apply. For
example, it is clear from the construction that if
$Q$ is surjective; respectively, injective, then so is $\b$, but
this follows also from general categorical principles: the
epimorphisms in $Vect$ are the surjective linear maps and in both
$Ban$ and
$Vect$, the monomorphisms are the injective morphisms.  On the
other hand, the epimorphisms in $Ban$ are not the surjective
operators but rather the operators with dense range and so need
not be epimorphisms in $Vect$; consequently, one would not expect
$\b$ to be an epimorphism in $Ban$  whenever $Q$ is (for an
example take $Q$ with dense proper range and $u$ so that $(uZ)\cap
(ZL)=\{0\}$--this forces $W$ to be $\{0\}$). 

From the canonical pullback construction it is also clear that if 
$Q$ has closed range, so does $\b$.  Thus if $Q$ is an
isomorphic embedding, so is $\b$.  The reason for taking the
$\ell_\infty$ sum of $L$ and $Z$ is that if $Q$ is an isometric
embedding  and $||u||\le 1$, then $\b$ is an isometric
embedding, and if
$Q$ is an isometric
quotient mapping and $||u||\le 1$, then $\b$ is an isometric
quotient mapping.  

 Notice also that
in (1.3) the kernels of $Q$ and $\b$ are isometrically
isomorphic, and in fact  (1.3)   can be extended to a
commutative diagram
$$\matrix
&0&\to&X&\to&L & {\buildrel Q\over \longrightarrow}  & Y \cr
&&&||& &\a \uparrow\phantom{\a} & & \phantom{u}\uparrow u
\cr 
&0&\to&V&\to& W &{\buildrel \b \over\longrightarrow} & Z
\cr
&&&&&\uparrow&&\uparrow
\cr
&&&&&W_0&{=}&Z_0
\cr
&&&&&\uparrow&&\uparrow\cr
&&&&&0&&0\cr
\endmatrix\leqno{(1.4)}
$$
with exact rows and columns, which of course cannot necessarily
be completed to short exact sequences.  However, if the top row of
(1.4) can be extended to a short exact sequence, so can the
second row--this is another way of saying that $\b$ is
surjective when $Q$ is surjective.  

Sometimes one can easily determine whether a commutative diagram (1.3)
is a pullback of (1.1).  For example, if 
$$\matrix
&0&\to&X&\to&L & {\buildrel Q\over \longrightarrow}  & Y &\to & 0\cr
&&&||& &\a \uparrow\phantom{\a} & & \phantom{u}\uparrow u
\cr 
&0&\to&V&\to& W &{\buildrel \b \over\longrightarrow} & Z &\to & 0
\cr
\endmatrix\leqno{(1.5)}
$$
is a  commutative exact
diagram in $Ban$, then (1.3)
is a pullback of (1.1).  Also, if  (1.4) is commutative and exact, and
$QL\cap uZ=Q\a W$--which is automatic when
$Q$ and $\b$ are surjective--then (1.3)  is a pullback of (1.1). 
Notice that it is enough to check these assertions in the nice category 
$Vect$, for then the unique bounded linear operator from $W$ to the
corner of the canonical pullback of (1.1) which makes the relevant
diagram commute  must be a surjective vector space isomorphism,
hence a surjective isomorphism in $Ban$ by the open mapping
theorem.  

Proposition~1.5, the first part of which is important for [KP],
says that the pullback construction provides an alternate way of
looking at the problem of when an operator factors or locally factors
through a quotient mapping:

\proclaim {Proposition 1.5} Consider the exact commutative diagram
$$\matrix
&0&\to&X&\to&L & {\buildrel Q\over \longrightarrow}  & Y &\to &0\cr
&&&||& &\a \uparrow\phantom{\a} & & \phantom{u}\uparrow u
\cr 
&0&\to&V&\to& W &{\buildrel \b \over\longrightarrow} & Z&\to &0
\cr
\endmatrix\leqno{(1.6)}
$$
\item{(1)} $u$ lifts to $L$ if and only if the second row splits.
\item{(2)} $u$ locally lifts to $L$ if and only if the second row
locally splits.
\endproclaim

\demo{Proof}  The ``if" direction is obvious both in (1) and (2).
Assume now that $u$  lifts to $L$, say  $u=Q\gam$, where $\gam:Z\to
L$.  From the discussion prior to the statement of Proposition~1.5,
we can assume without loss of generality that (1.3), the right
square of (1.6), is the canonical pullback of (1.1), the upper
right triangle of (1.6).  This makes it easy to define a lifting,
$\tau$, of $I_Z$ to $W$; namely, set $\tau z=(\gam z, z)$. Since
$\b$ is the projection onto the second component, this gives (1).

Part (2) follows from (1) by taking second adjoints in (1.6) and
applying Corollary~1.4(4) and Proposition~1.2(4). Alternatively, let
$w:C_1:\to Z$ be any operator and extend (1.6) via the pullback
construction to a commutative diagram with exact rows:

$$\matrix
&0&\to&X&\to&L & {\buildrel Q\over \longrightarrow}  & Y &\to &0\cr
&&&||& &\a \uparrow\phantom{\a} & & \phantom{u}\uparrow u
\cr 
&0&\to&V&\to& W &{\buildrel \b \over\longrightarrow} & Z&\to &0
\cr
&&&||& &  \uparrow  & & \phantom{w}\uparrow w
\cr
&0&\to&  V_1&\to&  W_1 &\longrightarrow  &
C_1&\to &0
\cr
\endmatrix\leqno{(1.7)}
$$

By hypothesis, the operator $uw$ factors through $Q$, so by part
(1) of Proposition~1.5  the bottom row of (1.7) splits,  hence $w$
factors through $\b$.  
This gives (2).\hfill\qed\enddemo

It is of considerable interest to determine when a locally
splitting short exact sequence must split. The
sequence
 $0\to X \to X^{**}\to X^{**}/X\to 0$, where $X \to X^{**}$ is the
canonical injection, must locally split--this immediate consequence of
Proposition~1.4 has long been known--and splits if
and only if $X$ is complemented in some dual space.  This can be
used to prove the following fact, which is a version of what is
called in [KP] Lindenstrauss' lifting criterion.

\proclaim {Lemma~1.6} Consider the diagram
$$\matrix
&0&\to&X&\to&L & {\buildrel Q\over \longrightarrow}  & Y &\to &0\cr
&&&&&&&\phantom{u}\uparrow u\cr
&&&&&&&Z\cr
\endmatrix
\leqno{(1.8)}
$$
where the top row is exact.  If $u$ locally lifts to $L$ and $X$ is
complemented in $X^{**}$, then $u$ lifts to $L$.
\endproclaim

\demo {Proof} Lindenstrauss' argument [Lin] provides a simple 
enough proof, but it is even easier to use Proposition~1.5.  Extend
(1.8) to (1.6).  The second row of (1.6) locally splits by
Proposition~1.5.  Now look at the commutative diagram

$$
\matrix
&0&\to&V^{**}&\to& W^{**}&{\buildrel \b^{**} \over\longrightarrow} &
Z^{**}&\to &0
\cr
&&&\uparrow& &  \uparrow  & &  \uparrow  
\cr
&0&\to&V&\to& W &{\buildrel \b \over\longrightarrow} & Z&\to &0
\cr
\endmatrix
\leqno{(1.9)}
$$
where the vertical arrows are the canonical embeddings and the rows
are exact.  The top row of (1.9) splits by Corollary~1.4.  The
space $V$, being isomorphic to $X$, is complemented in some dual
and hence in $V^{**}$, so the bottom row of (1.9) splits.  So $u$
factors through $Q$ by the trivial direction of Proposition~1.5 (2).
\hfill\qed\enddemo

Unfortunately, Lemma~1.6 is not of much use in determining when
a coextension of a $C(K)$ space must split; this is a problem
closely connected to the investigation of the so-called ``extension
property" considered in [JZ1], [JZ2].  

It was shown  in [KP] that pullbacks provide a quick proof (which,
however, relies on deep results from Banach space theory) of the
following result:

\proclaim {Proposition~1.7} If $Z^*$ has   cotype two and
$\ell_2$ is a quotient of $Z$; in particular, if $Z=C[0,1]$; then
there is an extension of $Z$ by $\ell_2$ which does not split.
\endproclaim

\demo {Proof} It is known [ELP], [KPec] that there is an extension,
$L$, of $\ell_2$ by $\ell_2$ which does not split.  Using the
pullback construct, we get a commuting diagram   
$$
\matrix
&0&\to&\ell_2&\to& L&{\buildrel Q \over\longrightarrow} &
\ell_2&\to &0
\cr
&&&||& &  \phantom{a}\uparrow \a  & & \phantom{u} \uparrow  u
\cr
&0&\to&\ell_2&\to& W &{\buildrel \b \over\longrightarrow} & Z&\to &0
\cr
\endmatrix
\leqno{(1.10)}
$$
with exact rows, $u$ a surjection, and the top row not splitting. 
The mapping $\a$ is a surjection since $u$ is, hence $\a^*$ is an
isomorphic embedding of $L^*$ into $W^*$.  If the bottom row splits,
then $W^*$ is isomorphic to the direct sum of $\ell_2$ and $Z^*$,
hence $W^*$, whence also $L^*$, has cotype two.  But for the known
constructions of such $L$'s, $L^*$ does not have cotype two. 
Actually, for any such $L$, 
$L^*$ cannot have cotype two since  that would force $L$ to have
type two (by Pisier's theorem [Pis1] and the Maurey--Pisier duality
theory [MP] for type--cotype in $K$-convex spaces),  in which case  
$0\to \ell_2\to L \to \ell_2\to 0$ would split by Maurey's
factorization theorem [Mau1]. \hfill\qed\enddemo   

Since every Banach space has quotients uniformly isomorphic to
$\ell_2^n$ for all $n$,  local Banach space theory
considerations show that the hypothesis in Proposition~1.7 that
$\ell_2$ be a quotient of $Z$ is not needed; it is that version
which appears in [KP]. In [KP] it is also noted that if  $Z$
contains subspaces uniformly isomorphic to $\ell_\infty^n$, then
there is an extension of $Z$ by
$\ell_2$ which does not split. This is also a consequence of 
Proposition~1.7 and local theory techniques.
 
We turn to the notion of pushout, which is dual (in the sense of
category theory) to that of pullback.  A commutative diagram

$$\matrix
 &\phantom{u}X & {\buildrel J\over \longrightarrow}  & L \cr
&u\downarrow &  &\phantom{\b} \downarrow \b \cr
&\phantom{u} Z &{\buildrel \a\over \longrightarrow}  &W \cr
\endmatrix\leqno{(1.11)}
$$
is called a {\sl pushout\/} of (1.2) provided that for every
commutative diagram
$$\matrix
 &\phantom{u}X & {\buildrel J\over \longrightarrow}  & L \cr
&u\downarrow &  &\phantom{\b} \downarrow \b_1 \cr
&\phantom{u} Z &{\buildrel \a_1\over \longrightarrow}  &W_1 \cr
\endmatrix\leqno{(1.12)}
$$
there is a unique morphism $W{\buildrel w\over \rightarrow} W_1$ so
that $\a_1=w\a$ and $\b_1=w\b$.  Pushouts are unique in an obvious
sense whenever they exist.  Pushouts exist in $Ban$.  Having
thought through pullbacks, one can build pushouts by first taking the
adjoint of (1.2) and constructing the pullback of it:
$$\matrix
 &\phantom{u}X^* & {\buildrel J^*\over \longleftarrow}  & L^* \cr
&u^*\uparrow \phantom{u^*}&  &\phantom{\b} \uparrow \tilde{\b}
\cr & Z^* &{\buildrel \tilde{\a} \over \longleftarrow} 
&\widetilde{W} 
\cr
\endmatrix\leqno{(1.13)}
$$
The subspace $\widetilde{W}=\{(z^*,x^*) : u^*z^*=J^*x^*\}$   is
weak$^*$ closed in $Z^*\oplus_\infty L^*=\left(Z\oplus_1 L\right)^*$
and the coordinate projections $\tilde{\a}$ and $\tilde{\b}$ are
weak$^*$ continuous.  Moreover, if one takes the adjoint of the
commutative diagram (1.12) and writes down the unique operator
$W_1^*\to
\widetilde{W}$ which makes the relevant diagram commute, one sees
that 
$W_1^*\to \widetilde{W}$ is weak$^*$ continuous.  Thus the preadjoint
of (1.13) is a indeed a pushout of (1.2), which we call the
{\sl canonical pushout;}  it is defined directly by setting
$W=\left(Z\oplus_1 L\right)/\overline{K}$, where 
$K=\{(ux,-Jx) : x\in X\}$, with $\a$ and
$\b$ the compositions of the natural mappings of $Z$ and $L$ into 
$Z\oplus_1 L$ with the quotient map from 
$Z\oplus_1 L$ onto $W$.  This is  the construction of the
pushout of (1.2) in $Vect$ only when $K$ is closed.  A natural
condition to guarantee that $K$ be closed is that $J$ have closed
range. Actually, the case where $J$ is even an isomorphic
embedding may be the only one considered in the Banach space
literature; at any rate, it is this case which has  played an
important role in Banach space theory.  The first deep application I
am aware of was due to Kisliakov [Kis].  The construction was also 
critical for Pisier's fundamental paper [Pis2].  Of course,
canonical pushouts play a major role in the  Kalton--Pe\l czy\'nski
paper [KP].  However, the categorical aspects of the canonical
pushout seem not to have been explicitly noted.   

Suppose that (1.11) is a pushout of (1.2).  Either directly from
the canonical construction or by taking adjoints and using the
pullback theory, one checks basic facts: If $J$ is surjective or an
isomorphic embedding or has closed range, then $\a$ has the same
property.  If
$||u||\le 1$ and $J$ is an isometric quotient map (respectively, an
isometric embedding), then $\a$ is an isometric quotient map (respectively, an
isometric embedding).  If $J$ is an epimorphism in $Ban$ (that is,
has dense range), so is $\a$.  The map $\a$ need not be injective
when $J$ is (take $J$ injective with dense proper range and let $u$
be a linear functional in $X^*$ which is not in $J^*L^*$).

The dual to the extension of the pullback diagram (1.4) (after
relabeling to agree with our pullback notation) is:
$$\matrix
 &\phantom{u}X & {\buildrel J\over \longrightarrow}  & L
&\longrightarrow& Y&\to &0\cr 
&u\downarrow &  &\phantom{\b} \downarrow
\b&&||
\cr &\phantom{u} Z &{\buildrel \a\over \longrightarrow}  &W 
&\longrightarrow &Y_1&\to &0\cr
\endmatrix\leqno{(1.14)}
$$
In order for the first row of (1.14) to be exact, $J$ must have
closed range, in which case if (1.11)
is a pullback diagram $\a$ also has closed range and the quotients 
$Y\equiv L/JX$ and 
$Y_1\equiv W / \a Z$  are naturally isomorphic.  

If (1.4) is a pullback diagram
and both $J$ and $u$ have closed range, then the pushout
construction produces the dual commutative diagram to (1.5)
$$\matrix
 &\phantom{u}X & {\buildrel J\over \longrightarrow}  & L
&\longrightarrow &Y&\to &0\cr 
&u\downarrow \phantom{u} &  &\phantom{\b} \downarrow
\b&&||
\cr & Z &{\buildrel \a\over \longrightarrow}  &W 
&\longrightarrow& Y_1&\to& 0
\cr&\downarrow&&\downarrow\cr
&V&=&V_1\cr
&\downarrow&&\downarrow\cr
&0&&0\cr
\endmatrix\leqno{(1.15)}
$$
where the columns and rows are exact.  If $J$ is an isomorphic 
embedding, the
rows in (1.14) can be extended to short exact sequences.  
When both $J$ and
$u$ are isomorphic embeddings, we get the commuting diagram
 $$\matrix
&&&0&&0\cr
&&&\downarrow&&\downarrow\cr
&0&\to &X & {\buildrel J\over \longrightarrow}  & L
&\longrightarrow& Y&\to &0\cr 
&&&u\downarrow \phantom{u} &  &\phantom{\b} \downarrow
\b&&||
\cr
&0&\to & Z &{\buildrel \a\over \longrightarrow}  &W 
&\longrightarrow &Y_1&\to& 0
\cr
&&&\downarrow&&\downarrow\cr&&
&V&=&V_1\cr&&
&\downarrow&&\downarrow\cr&&
&0&&0\cr
\endmatrix\leqno{(1.16)}
$$
where the rows and columns are exact; this is used in section 2.

Also used in section 2 is part (2) of Proposition~1.8 (part (1) is
important for [Kis], [Pis2], and [KP]), which is dual to
Proposition~1.5.

\proclaim{Proposition 1.8} Consider the commutative exact diagram
 $$\matrix
&0&\to &\phantom{u}X & {\buildrel J\over \longrightarrow}  & L
&\longrightarrow Y\to 0\cr 
&&&u\downarrow \phantom{u} &  &\phantom{\b} \downarrow
\b&||
\cr
&0&\to & Z &{\buildrel \a\over \longrightarrow}  &W 
&\longrightarrow Y_1\to 0
\cr
\endmatrix
\leqno{(1.17)}
$$
\item{(1)} $u$ extends to $L$ if and only if the second row splits.
\item{(2)} $u$ colocally extends to $L$ if and only if the second
row locally splits.
\endproclaim

\demo{Proof of  (2)}If $u$ colocally extends to $L$, then by
Proposition~1.1 $u^*$ factors through $J^*$, so that by
Proposition~1.5 (take the adjoint diagram of (1.16)), \break
$0\leftarrow Z^*{\buildrel \a^*\over\longleftarrow}W^*\leftarrow
Y^*_1\leftarrow 0$ splits.  Hence also $0\to Z^{**} {\buildrel
\a^{**}\over \longrightarrow W^{**}} \to Y_1^{**} \to 0$ splits,
whence again by Corollary~1.4, $0\to Z {\buildrel
\a \over \longrightarrow} W  \to Y_1  \to 0$ locally splits.

Conversely, if $0\to Z {\buildrel
\a \over \longrightarrow} W  \to Y_1  \to 0$ locally splits, then
the sequence 
\break
$0\leftarrow Z^*{\buildrel \a^*\over\longleftarrow}W^*\leftarrow
Y^*_1\leftarrow  0$ splits by Corollary~1.4.  Hence by
Proposition~1.5,
$u^*$ factors through $J^*$, whence $u$ colocally extends to $L$ by
Proposition~1.1. \hfill\qed\enddemo

\vfill\eject
 
\heading 2. Extensions of $c_{0}$. 
\endheading

In this section we prove:

\proclaim {Theorem 2.1}  Suppose that 
$$0\to X \to L \to Y \to 0\leqno{(2.1)}$$ is
exact with $X$ separable, $c_0$ is isomorphic to a subspace of
$Y$, and
$L$ embeds into  a Banach lattice which does not contain
$\ell_\infty^n$'s uniformly. Suppose
$$0\to X \to Z \to V \to 0\leqno{(2.2)}$$ 
is a locally splitting short exact sequence with $Z$  separable 
and $Z$ embeds into a Banach lattice which does not contain 
$\ell_\infty^n$'s uniformly. Then
$Z$ is not complemented in its bidual.\endproclaim

This theorem has a corollary which can
be stated in the language of Banach space theory as:

\proclaim {Corollary 2.2} Suppose that $L$  embeds into a Banach
lattice which does not contain
$\ell_\infty^n$'s uniformly and $Q$ is an operator from $L$ onto 
some Banach space $Y$.   If 
$\ker Q$ has GL-{\sl l.u.st.,}
 then $Y$ does not contain 
$\ell_\infty^n$'s uniformly.
 \endproclaim

If the conclusion in Corollary~2.2 is weakened to ``$c_0$ does
not embed into $Y$", the resulting statement (at least when $X$
is separable) is immediate from Theorem~2.1 and known results. In
the appendix we show how to deduce Corollary~2.2 from Theorem~2.1,
preferring in this section to concentrate on the proof of
Theorem~2.1 itself.   The most important case is $L=L_1$, but this
case is  not easier than the general one. However, the case
$L=L_1$ does lend an easier proof that $X$ does not have {\sl l.u.st.}
in its original sense, and we mention in the proof of Theorem~2.1
how to streamline the proof to obtain just this.

Note that even when $Y=c_0$, the hypotheses on $L$ in Theorem~2.1
cannot be replaced by the  conditions that that $c_0$ does not embed
into
$L$ and $L$ is itself a lattice.  Indeed, it is clear that
the identity on $c_0$ locally factors through the natural quotient 
map from $\left(\sum_{n=1}^\infty \ell_\infty^n\right)_1$ onto $c_0$.

\demo {Proof of Theorem~2.1} We can assume that $X$ is 
a subspace of $L$ and that
$Y=L/X$. It is easy to see that there is
a separable superspace $L'$ of $X$ in
$L$ so that $L'/X$ is isomorphic to $c_0$. So by replacing $L$
by $L'$, we can assume  that
$L$, and {\sl a fortiori\/} also $Y$,  are separable.  We are thus
considering the  exact diagram:
  
 $$\matrix
&&&0\cr&&&\downarrow\cr
&0&\to &X & {\buildrel J\over \longrightarrow}  & L
&\longrightarrow &Y&\to& 0\cr 
&&&u\downarrow \phantom{u} 
\cr
&& & Z
\cr
&&&\downarrow\cr&&&V\cr
&&&\downarrow\cr&&
&0\cr
\endmatrix\leqno{(2.3)}
$$
where $L$ and $Z$ both embed into  separable Banach lattices
which do not contain 
$\ell_\infty^n$'s uniformly, $Y$ contains a  copy of
$c_0$ (which, since $Y$ is separable, is necessarily complemented
because $c_0$ has the separable extension property ([Sob], 
[LT1,~Th.~2.f.5]), and the column locally splits.  In some sense,
the main point is one that is obvious to any self-respecting
algebraist:  {\sl   to study\/}
\rm (2.3), {\sl one should complete it to a full diagram.}  But
since
$J$ and $u$ are both isomorphic embeddings, the pushout
construction extends (2.3) to the exact commutative diagram
(1.16), which we repeat as (2.4):

 $$\matrix
&&&0&&0\cr
&&&\downarrow&&\downarrow\cr
&0&\to &X & {\buildrel J\over \longrightarrow}  & L
&{\buildrel q \over \longrightarrow}& Y&\to &0\cr 
&&&u\downarrow \phantom{u} &  &\phantom{\b} \downarrow
\b&&||
\cr
&0&\to & Z &{\buildrel \a\over \longrightarrow}  &W 
&{\buildrel q_1\over \longrightarrow}& Y_1&\to& 0
\cr
&&&\downarrow&&\downarrow\cr&&
&V&=&V_1\cr&&
&\downarrow&&\downarrow\cr&&
&0&&0\cr
\endmatrix\leqno{(2.4)}
$$

Before proceeding further, let us see that when $L=L_1$, the
space $X$ cannot have {\sl l.u.st.}; this answers the question
Kalton and Pe\l czy\'nski posed in their lectures on [KP].  If
$X$ has {\sl l.u.st.}, then it is known that there is a locally
splitting short exact sequence (2.2) with
$Z$ a separable Banach lattice which isomorphically embeds into
$L_1$ (see the appendix). In view of   Corollary~1.4, this
implies  that
$V^{**}$ embeds into $Z^{**}$ which embeds into the abstract
$L_1$-space $L_1^{**}$, and hence by [LP] $V$ embeds into
$L_1$.  Since the first column in (2.4) locally splits, so does
the second. But then $W^{**}$ is isomorphic to $L_1^{**}\oplus
V^{**}$ and hence $W$ also embeds into $L_1$.  Now we need a key
 analytical lemma proved, but not stated, in [KP]. (See the proof
of Proposition~2.2 in [KP]. Actually, in [KP]  $P$ is constructed
so that $PQ$ is  even representable; that is, factors through
$\ell_1$.)  Later we prove a generalization of  Lemma~2.3 
in order to prove Theorem~2.1.

\proclaim{Lemma 2.3} If $W$ is a subspace of $L_1$ and $q_1$ is an
operator from $W$ into a space $Y_1$ which contains a complemented
copy of $c_0$, then there is a projection $P$ on $Y_1$ with
$PY_1$ isomorphic to $c_0$ such that $Pq_1$ is {\sl completely
continuous;\/} that is, carries weakly convergent sequences into
norm convergent sequences.
\endproclaim

Let $P$ be given from Lemma~2.3 and let $v$ be an operator from
$L_1$ into the isomorph $PY_1 $ of $c_0$ which is not completely
continuous; there are a wealth of such operators.  $v$ cannot
factor through the completely continuous operator $Pq_1$, which is
to say that $v$
 cannot factor through $q_1$.     But $v$ locally factors through
$q_1$ since $q_1$ is surjective.  By Lemma~1.6, $Z$ is not
complemented in its bidual.  But every Banach lattice which does
not contain an isomorph of $c_0$ is complemented in its bidual
via a band projection [LT2].  

This digression from the proof of Theorem~2.1 in fact motivates
the proof, to which we now return.  Since the properties of being
a Banach lattice and of not containing $\ell_\infty^n$'s
uniformly are both preserved under passage to biduals, just as in
the digression we conclude that $W$ embeds into a Banach lattice
which does not contain $\ell_\infty^n$'s uniformly.  The further
argument in the digression shows that in order to complete the
proof of Theorem~2.1, it is enough to verify the following
generalization of Lemma~2.3:

\proclaim{Lemma 2.4} If $W$ is a separable subspace of a Banach
lattice which does not contain $\ell_\infty^n$'s uniformly and
$q_1$ is an operator from $W$ into a space $Y_1$ which contains a
complemented copy of $c_0$, then there is a projection $P$ on
$Y_1$ with $PY_1$ isomorphic to $c_0$ such that $Pq_1w$ is   
completely continuous for every operator $w$ from $L_1$ into $W$.
\endproclaim

For the proof of Lemma~2.4 we need one sublemma and a couple of
known facts.

\proclaim{Sublemma 2.5} Let $W$ be a separable Banach
lattice which does not contain $\ell_\infty^n$'s uniformly. 
Suppose that $\fn$ is a weak$^*$ null sequence in $W^*$.  Then
there exist $g_n$ in the convex hull of $\{f_k\}_{k=n}^\infty$ so
that $\absgn$ is weak$^*$ null in $W^*$.
\endproclaim

\demo{Proof} It is known that there exists $q<\infty$, a measure $\mu$,
and an operator $u$ from $L_q(\mu)$ into $W$ with dense range which is an
interval preserving lattice homomorphism.  Since standard texts do not
 include this fact, here is a sketch of the proof (unexplained terminology
as well as the quoted theorems about lattices can be found in [LT2]): 
Let
$x$ be a weak order unit for $W$ and let $X$ be the linear span of the order
interval
$[-x,x]$ with  $[-x,x]$ as its unit ball.  $X$ is then an abstract
$M$-space and so can be identified, as a Banach lattice, with $C(K)$ for
some compact Hausdorff space $K$ by Kakutani's representation theorem. 
Since $W$ does not contain
$\ell_\infty^n$'s uniformly, the injection $j$ from $X$ into $W$ is
$q$-summing for some $q$ by a theorem of Maurey's [DJT,~p.~223].  Choose
a Pietsch measure [DJT,~p.~45] $\mu$ for $u$; then $j$ factors through
the natural injection $i$ from $X=C(K)$ into $L_q(\mu)$; say, $j=ui$
where 
$L_q(\mu){\buildrel u\over \longrightarrow} W$.  The operator $u$ is of
course uniquely defined because $C(K)$ is dense in  $L_q(\mu)$; that $u$
has the stated properties can be deduced from the fact that $j$ has those
properties.

Taking adjoints, we see that $W^* {\buildrel u^*\over \longrightarrow}
L_p(\mu)$  (${1\over p}+{1\over q}=1$) is an injective (since $u$ has dense
range) lattice homomorphism (since $u$ is interval preserving; see
[AB,~p.~92]).  By the weak$^*$ continuity of $u^*$,  $u^*f_n \to 0$ weakly
in   
$L_p(\mu)$, and hence there exist $g_n$ in the convex hull of
$\{f_k\}_{k=n}^\infty$ so
that $||u^*g_n||_p \to 0$ and hence $||\,|u^*g_n|\,||_p\to 0$.  But 
$|u^*g_n| = u^*|g_n|$ because $u^*$ is a lattice homorphism.  But then
the only possible weak$^*$ cluster point in $W^*$ of $\absgn$ is
$0$, so that in fact $\absgn$ must converge weak$^*$ in $W^*$ to $0$.
\enddemo

\proclaim{Fact 2.6} If $W$ is a Banach lattice which does not
contain a subspace isomorphic to $c_0$ and $w$ is an operator
from $L_1$ into $W$, then $w\Ball(L_\infty)$ is order bounded.
\endproclaim
In fact, the stated hypothesis implies that every operator from 
$L_1$ into $W$ has a modulus [AB,~p.~249].

\proclaim{Fact 2.7} An operator $w$ from $L_1$ is completely
continuous if and only if $w\Ball(L_\infty)$ is relatively
compact. \endproclaim

Fact~2.7 can be found in [Uhl].  
Actually, we do not really need it because it is evident
that there are operators from $L_1$ into $c_0$ which take the
ball of $L_\infty$ into a non-relatively compact set.

\demo{Proof of Lemma 2.4} It is clear that there is no loss of
generality in assuming that $Y_1=c_0$.   
  By replacing $W$ by the closure of the
(necessarily separable) sublattice it generates in some containing
Banach lattice which does not contain $\ell_\infty^n$'s
uniformly, we can assume, since $c_0$ has the separable
extension property, that
$W$ itself is a Banach lattice.

Set $f_n=q_1^*e_n^*$, where
$\{e_n^*\}_{n=1}^\infty$ is the unit vector basis of $\ell_1=c_0^*$.
Clearly $\fn$ tends weak$^*$ in $W^*$ to $0$, so
from Sublemma 2.5 we get  
$n_1<n_2<n_3 \dots$ and  
$d_n^*\in
\hbox{  co\,} \{e^*_i\}_{i={k_n}+1}^{k_{n+1}}$ so that, setting
$g_n=q_1^*d_n^*$, 
$ \absgn$ converges weak$^*$ in $W^*$ to $0$.  Define 
$x_n=\sum_{i=k_{n}+1}^{k_{n+1}} e_i$ in $c_0$ and let $P$ be the
contractive projection on  $c_0$ defined by $P=\sum d^*_n\otimes
x_n$.  

Let $w$ be any operator from $L_1$ into $W$.  By Fact 2.7, in
order to check that $Pq_1w$ is completely continuous, it is
enough to show that $Pq_1w\Ball (L_\infty)$ is relatively compact.
This amounts to checking that 
$$
\sup_{g\in w Ball (L_\infty)} |\langle g, g_n\rangle| \to 0
{\hbox{ \ as \ } } n\to \infty.\leqno{(2.5)}
$$
But by Fact~2.6, there
exists $h\ge 0$ in $W$ so that $w\Ball (L_\infty)$ is contained
in the order interval $[-h,h]$.  We thus have for each $n$:
$$
\sup_{g\in w Ball (L_\infty)} |\langle g, g_n\rangle|\le 
\langle |h|, |g_n|\rangle.
$$
Since $\absgn$ is weak$^*$ null in $W^*$, (2.5) follows.
This completes the proof of Lemma~2.4 and hence also the proof of
Theorem~2.1. \hfill\qed\enddemo\enddemo

\noindent {\bf Remark  2.8} If the Banach space $X$ fails to have
GL-{\sl l.u.st.} and $X\subset Y \subset X^{**}$, then the identity on
$X$  colocally extends to  $Y$ by Corollary 1.4 (2) and hence $Y$ also
fails GL-{\sl l.u.st.}. 

 If $S$ is a Sidon subset of the compact
Abelian group $G$, then $L_1(G)/L^1_{\tilde{S}}(G)$ is isomorphic to
$c_0$ (see [KP] for background). 
 Pe\l czy\'nski pointed out that if we apply 
Remark 2.8 and Corollary 2.2 to this kind of example,  we obtain the
following information about the classical object $M_T(G)$,  the set of
finite measures on $G$ whose Fourier transforms are supported on
$T$:

\proclaim {Corollary 2.9} If $S$ is a Sidon subset of the compact
Abelian group $G$, then $M_{\tilde{S}}(G)$ does not have GL-{\sl
l.u.st.}.
\endproclaim

\vfill\eject

\heading 3. Appendix. 
\endheading

The main background needed for deriving Corollary~2.2 from Theorem~2.1 is
the following theorem, which is  a restatement of results from [FJT]:

\proclaim {Theorem 3.1}  A Banach space $X$ has GL-{\sl l.u.st.} if and
only if there is a locally splitting short exact sequence 
$0\to X \to Z \to V \to 0$ with $Z$ a Banach lattice.  Moreover, if in
addition to having GL-{\sl l.u.st.}
$X$ does not contain $\ell_\infty^n$'s uniformly, then $Z$ may be chosen
not to contain $\ell_\infty^n$'s uniformly.  Also,  if $X$ has {\sl
l.u.st.},
$Z$ may be chosen to be finitely crudely representable in $X$ (that is,
the finite dimensional subspaces of $Z$ embed into $X$ with uniformly
bounded isomorphism constants).
\endproclaim

In [FJT] and also [Mau2] it was remarked that a space $X$ has 
GL-{\sl l.u.st.} if and only if $X^{**}$ is complemented in a Banach
lattice (see [DJT,~p.~348] for a proof, but keep in mind that in [DJT] 
GL-{\sl l.u.st.} is called {\sl l.u.st.} while {\sl l.u.st.} is called 
DPR-{\sl l.u.st.}).  This gives the first statement in Theorem~3.1.  The
``also" statement is a consequence of Corollary~2.2 in [FJT], while
the  ``moreover" follows from Proposition~2.6(i) and Remark~2.8 of [FJT]. 
Of course, for us the definition of GL-{\sl l.u.st.} is irrelevant,
since in section~2 we use only the characterization given by Theorem~3.1.

Notice that the ``also" statement yields that if $X$ is a subspace
of an $L_1$ space and $X$ has {\sl l.u.st.}, then the Banach lattice $Z$
from Theorem~3.1 can be taken finitely crudely representable in $L_1$,
hence $Z$ embeds into an $L_1$ space by [LP].  This was used in the
``digression" part of the proof of Theorem~2.1.

\demo {Proof of Corollary~2.2} The reader who is not familiar with
 ultrapowers of Banach spaces will find enough in chapter 8 of [DJT] to
make the verification of claims we make about ultrapowers easy. The
statements we make about GL-{\sl l.u.st.} are probably more obvious from the
definition than from the equivalent form given by the first statement in
Theorem~3.1; again, [DJT] is sufficient reference.

The property of $L$ in Corollary~2.2;
namely, that $L$ embeds into a Banach lattice which does not contain 
$\ell_\infty^n$'s uniformly; is stable under the taking of ultrapowers,
as is  the property of having GL-{\sl l.u.st.}.  On the other hand, any
ultrapower of a  space which contains $\ell_\infty^n$'s uniformly must
contain $c_0$.  Consequently, to prove Corollary~2.2 it is enough to check
that $c_0$ does not embed into $Y$.  Assume, for contradiction, that $c_0$ 
does embed into $Y$. Choose a separable subspace $L'$ of $L$ so that
$QL'$ is closed and contains a copy of $c_0$.  Since $\ker Q$ has
GL-{\sl l.u.st.}, there is a separable subspace $X$ of $\ker Q$ containing
the intersection of $\ker Q$ with $L'$ and which has GL-{\sl l.u.st.}.  By
replacing $L'$ with the closed span of $L'\cup X$, we can assume that 
$X=L'\cap \ker Q$.  Thus we have a short exact sequence
$0\to X \to L' \to Y_0 \to 0$  with $L'$ a separable subspace of some Banach
lattice which does not contain $\ell_\infty^n$'s uniformly, $X$ has
GL-{\sl l.u.st.}, and $Y_0$ contains a copy of $c_0$.  But by Theorem~3.1,
there is a locally splitting short exact sequence 
$0\to X \to Z \to V \to 0$ with $Z$  a Banach lattice which does not
contain $\ell_\infty^n$'s uniformly.  Moreover, by replacing $Z$ with the
 closed sublattice generated by $X$, we can assume that $Z$ is
separable.  The lattice $Z$ is complemented in its bidual because it does
not contain a copy of $c_0$ [LT2].  This contradicts Theorem~2.1 and
completes the proof of Corollary~2.2.\hfill\qed\enddemo

\vfill\eject

\enddocument


\Refs
\widestnumber\key{KPec}
\def\n{\key}



\ref\n{AB}
\by  C. D. Aliprantis and O. Burkinshaw
\book  Positive operators
\publ  Academic Press 
\yr 1985 
\endref

\ref\n{Bou}
\by  J. Bourgain
\paper A counterexample to a complementation  problem
\jour Compositio Math.
\vol 43
\yr  1981
\pages 133--144
\endref

\ref\n{CG}
\by  J. M. F. Castillo and M. Gonz\'alez
\book  Three-space problems in Banach space theory 
\yr  1996
\endref


\ref\n{DJT}
\by  J. Diestel, H. Jarchow, and A. Tonge
\book  Absolutely summing operators
\publ  Cambridge University Press 
\yr  1995
\endref




\ref\n{Dom1}
\by  P. Doma\'nski
\book Extensions and liftings of linear operators
\publ  Adam Mickiewicz Univ.  
\publaddr Pozna\'n
\yr  1987
\endref

\ref\n{Dom2}
\by   P. Doma\'nski
\paper  On the splitting of twisted sums, and the three
space problem for local convexity
\jour Studia Math.
\vol 82 
\yr  1985
\pages 155--189
\endref

\ref\n{DPR}
\by  E. Dubinsky, A. Pe\l czy\'nski, and H. P. Rosenthal
\paper On Banach spaces $X$ for which 
$\Pi\sb{2}({\cal L}\sb{\infty}, X)=B({\cal L}\sb{\infty}, X)$
\jour Studia Math.
\vol 44
\yr  1972
\pages 617--648
\endref

\ref\n{ELP}
\by  P. Enflo, J. Lindenstrauss, and G. Pisier
\paper On the three space problem
\jour Math. Scand.
\vol 36
\yr  1975
\pages 199--210
\endref

\ref\n{FJT}
\by  T. Figiel, W. B. Johnson, and L. Tzafriri
\paper On Banach
lattices and spaces having local unconditional structure, with applications
to Lorentz function spaces
\jour J. Approximation Theory
\vol 13
\yr  1975
\pages 395--412
\endref


\ref\n{GL}
\by  Y.~Gordon and D.~R.~Lewis
\paper Absolutely summing operators 
and local  unconditional structures
\jour Acta Math.
\vol 133
\yr  1974
\pages 27--48
\endref


\ref\n{J1}
\by    W. B. Johnson 
\paper Factoring compact
operators
\jour Israel J. Math. 
\vol 9
\yr  1971
\pages 337--345
\endref

\ref\n{J2}
\by    W. B. Johnson 
\paper A complementably universal conjugate Banach space and its
relation to the approximation problem
\jour Israel J. Math. 
\vol  13
\yr   1972
\pages  201--310
\endref

\ref\n{JRZ}
\by    W. B. Johnson, H. P. Rosenthal, and M. Zippin 
\paper On bases, finite dimensional decompositions, and weaker
structures in Banach spaces 
\jour Israel J. Math. 
\vol  9
\yr   1971
\pages  488--506
\endref


\ref\n{JZ1}
\by    W. B. Johnson and M. Zippin
\paper Extension of operators from
subspaces of  $c_0(\gamma)$ into $C(K)$ spaces
\jour Proc.\ AMS 
\vol  107
\yr   1989
\pages  751--754
\endref



\ref\n{JZ2}
\by    W. B. Johnson and M. Zippin
\paper Extension of operators from weak$^*$-closed subspaces of
$\ell_1$
\jour Studia Math. 
\vol   117
\yr   1995
\pages   43--55
\endref

\ref\n{Kal}
\by  N. J.  Kalton 
\paper  Locally complemented subspaces and ${{\cal L}_{p}}$-spaces for
$0<p<1$
\jour Math. Nachr.
\vol  115
\yr   1984
\pages   71--97
\endref


\ref\n{KPec}
\by  N. J.  Kalton and N. T. Peck
\paper  Twisted sums of sequence spaces and the three space problem
\jour Trans. AMS
\vol  255
\yr   1979
\pages   1--30
\endref



\ref\n{KP}
\by  N. J.  Kalton and A. Pe\l czy\'nski
\paper  Kernels of surjections from ${\Cal L}_1$-spaces with an application
to Sidon sets
\jour Math. Annalen
\vol  
\yr   
\pages   
\endref

\ref\n{Kis}
\by  S. V. Kisliakov
\paper  On spaces with ``small" annihilators
\jour Sem. Leningrad Otdel. Math. Institute Steklov (LOMI)
\vol  65
\yr   1976
\pages   192--195
\endref


\ref\n{Lin}
\by    J. Lindenstrauss 
\paper On a certain subspace of $\ell_1$
\jour Bull. Acad. Polon. Sci. Ser. Sci. Math. Astronom. Phys. 
\vol  12
\yr   1964
\pages  539--542
\endref

\ref\n{LP}
\by    J. Lindenstrauss  and A. Pe\l czy\'nski
\paper  Absolutely summing
operators in ${\cal L}_p$ spaces and their applications
\jour Studia Math.
\vol  29
\yr   1968
\pages  275--326
\endref

\ref\n{LR}
\by  J. Lindenstrauss and H. P. Rosenthal
\paper  The ${\Cal L}_p$-spaces
\jour Israel J. Math.
\vol  7
\yr   1969
\pages  325--349
\endref

\ref\n{LT1}
\by  J. Lindenstrauss and L. Tzafriri
\book  Classical Banach spaces I, Sequence spaces
\publ  Sprin\-ger-Verlag
\publaddr Berlin
\yr  1977
\endref

\ref\n{LT2}
\by  J. Lindenstrauss and L. Tzafriri
\book  Classical Banach spaces II, Function spaces
\publ  Sprin\-ger-Verlag
\publaddr Berlin
\yr  1979
\endref

\ref\n{Mau1}
\by    B. Maurey
\paper Un th\'eor\'eme de prolongment
\jour C. R. Acad.  
Paris 
\vol  279
\yr   1974
\pages  329--332
\endref

\ref\n{Mau2}
\by    B. Maurey
\paper Type et cotype dans les espaces munis de structures locales
inconditionelles
\jour S\'eminaire Maurey-Schwartz 1973/74
\vol  
\yr   1974
\pages  Exp. 24--25
\endref


\ref\n{MP}
\by    B. Maurey and G. Pisier
\paper S\'eries de variables  al\'eatoires 
ind\'ependantes et propri\'et\'es g\'eom\'e\-triques des  espaces de
Banach
\jour Studia Math.   
\vol  58
\yr   1976
\pages  45--90
\endref



\ref\n{Pis1}
\by     G. Pisier
\paper Holomorphic semi-groups and the geometry of  
Banach spaces  
\jour Ann. of Math.   
\vol  115
\yr   1982
\pages  375--392
\endref




\ref\n{Pis2}
\by     G. Pisier
\paper Counterexamples to a conjecture of Grothendieck
\jour Acta Math.   
\vol  151
\yr   1983
\pages  181--209
\endref




\ref\n{Sob}
\by     A. Sobczyk
\paper Projections of the space $m$ onto its subspace $c_0$
\jour    Bull. AMS
\vol  47
\yr   1941
\pages  938--947
\endref

\ref\n{Uhl}
\by     J. J. Uhl
\paper Three quickies in the theory of vector measures
\jour    The Altgeld Book, 1975/76, \ The  University of Illinois 
Functional Analysis seminar
\yr   1976
\pages  
\endref

\endRefs
\bye